\documentclass[12pt,thmsa]{article}%
\usepackage{amssymb}
\usepackage{amsfonts}
\usepackage{sw20bams}%
\usepackage{amsmath}%
\usepackage{graphicx}
\providecommand{\U}[1]{\protect\rule{.1in}{.1in}}
\begin{document}

\author{Steven R. Finch}
\title{Mean Width of a Regular Simplex}
\date{March 12, 2016}
\maketitle

\begin{abstract}
The mean width is a measure on $n$-dimensional convex bodies. An integral
formula for the mean width of a regular $n$-simplex appeared in the electrical
engineering literature in 1997. As a consequence, expressions for the expected
range of a sample of $n+1$ normally distributed variables, for $n\leq6$, carry
over to widths of regular $n$-simplices. As another consequence, precise
asymptotics for the mean width become available as $n\rightarrow\infty$.

\end{abstract}

\footnotetext{Copyright \copyright \ 2011, 2016 by Steven R. Finch. All rights
reserved.}Let $C$ be a convex body in $\mathbb{R}^{n}$. A\ \textbf{width} is
the distance between a pair of parallel $C$-supporting planes (linear
varieties of dimension $n-1$). Every unit vector $u\in\mathbb{R}^{n}$
determines a unique such pair of planes orthogonal to $u$ and hence a width
$w(u)$. Let $u $ be uniformly distributed on the unit sphere $S^{n-1}\subset$
$\mathbb{R}^{n}$. Then $w$ is a random variable and
\[
\mathbb{E}\left(  w_{3}\right)  =\frac{3}{2\pi}\arccos\left(  -\frac{1}%
{3}\right)
\]
for $C=$ the regular $3$-simplex (tetrahedron) in $\mathbb{R}^{3}$ with edges
of unit length and
\[
\mathbb{E}\left(  w_{4}\right)  =\frac{10}{3\pi^{2}}\left[  3\arccos\left(
-\frac{1}{3}\right)  -\pi\right]
\]
for $C=$ the regular $4$-simplex in $\mathbb{R}^{4}$ with edges of unit length
\cite{Bo, BS}. \ Our contribution is to extend the preceding \textbf{mean
width} results to regular $n$-simplices in $\mathbb{R}^{n}$ for $n\leq6 $. We
similarly extend the following \textbf{mean square width} result:%
\[
\mathbb{E}\left(  w_{3}^{2}\right)  =\frac{1}{3}\left(  1+\frac{3+\sqrt{3}%
}{\pi}\right)
\]
which, as far as is known, first appeared in \cite{Fi1}.

The key observation underlying our work is due to Sun \cite{Sn}, which in turn
draws upon material in \cite{Wb, Bk}. \ It does not seem to have been
acknowledged in the mathematics literature. \ After most of this paper was
written, we found \cite{HZ}, which assigns priority to to Hadwiger \cite{Hd}
and to Ruben \cite{Rb} for closely related ideas.

\section{Order Statistics}

Let $X_{1},$ $X_{2},$ $...,$ $X_{n}$ denote a random sample from a Normal
$(0,1)$ distribution, that is, with density function $f$ and cumulative
distribution $F$:
\[%
\begin{array}
[c]{ccc}%
f(x)=\dfrac{1}{\sqrt{2\pi}}\exp\left(  -\dfrac{x^{2}}{2}\right)  , &  & F(x)=%
{\displaystyle\int\limits_{-\infty}^{x}}
f(\xi)\,d\xi=\dfrac{1}{2}\operatorname*{erf}\left(  \dfrac{x}{\sqrt{2}%
}\right)  +\dfrac{1}{2}.
\end{array}
\]
The first two moments of the \textbf{range}
\[
r_{n}=\max\{X_{1},X_{2},...,X_{n}\}-\min\{X_{1},X_{2},...,X_{n}\}
\]
are given by \cite{Tp, Mr}%
\[
\mu_{n}=\mathbb{E}(r_{n})=%
{\displaystyle\int\limits_{-\infty}^{\infty}}
\left\{  1-F(x)^{n}-[1-F(x)]^{n}\right\}  dx,
\]%
\[
\nu_{n}=\mathbb{E}(r_{n}^{2})=2%
{\displaystyle\int\limits_{-\infty}^{\infty}}
{\displaystyle\int\limits_{-\infty}^{y}}
\left\{  1-F(y)^{n}-[1-F(x)]^{n}+[F(y)-F(x)]^{n}\right\}  dx\,dy.
\]
For small $n$, exact expressions are possible \cite{Da, Ru, Wa}:
\[%
\begin{array}
[c]{lll}%
\mu_{2}=\frac{2}{\sqrt{\pi}}=1.128..., &  &
\begin{array}
[c]{c}%
\nu_{2}=2,
\end{array}
\\
\mu_{3}=\frac{3}{\sqrt{\pi}}=1.692..., &  &
\begin{array}
[c]{c}%
\nu_{3}=2\left(  1+\frac{3\sqrt{3}}{2\pi}\right)  =3.653...,
\end{array}
\\
\mu_{4}=\frac{6}{\sqrt{\pi}}\left(  1-2S_{2}\right)  =2.058..., &  &
\begin{array}
[c]{c}%
\nu_{4}=2\left(  1+\frac{3+\sqrt{3}}{\pi}\right)  =5.012...,
\end{array}
\\
\mu_{5}=\frac{10}{\sqrt{\pi}}\left(  1-3S_{2}\right)  =2.325..., &  &
\begin{array}
[c]{c}%
\nu_{5}=2\left(  1+\frac{5\sqrt{3}}{2\pi}+\frac{30}{\pi}S_{1/2}-\frac
{5\sqrt{3}}{\pi}S_{3}\right)  =6.156...,
\end{array}
\\
\mu_{6}=\frac{15}{\sqrt{\pi}}\left(  1-4S_{2}+2T_{2}\right)  =2.534..., &  &
\begin{array}
[c]{c}%
\nu_{6}=2\left(  1+\frac{5(9+2\sqrt{3})}{2\pi}-\frac{90}{\pi}S_{2}%
-\frac{15\sqrt{3}}{\pi}S_{3}\right)  =7.142...,
\end{array}
\\
\mu_{7}=\frac{21}{\sqrt{\pi}}\left(  1-5S_{2}+5T_{2}\right)  =2.704..., &  &
\begin{array}
[c]{l}%
\nu_{7}=2\left(  1+\frac{35\sqrt{3}}{4\pi}+\frac{210}{\pi}S_{1/2}-\frac
{105}{\pi}S_{2}-\frac{35\sqrt{3}}{\pi}S_{3}\right. \\
\;\;\;\;\;\;\;\;\;\;\;\;\;\;\;\left.  +\frac{35\sqrt{3}}{2\pi}T_{3}+\frac
{210}{\pi}U-\frac{420}{\pi}V\right)  =8.007...
\end{array}
\end{array}
\]
where
\[
S_{k}=\frac{\sqrt{k}}{\pi}%
{\displaystyle\int\limits_{0}^{\pi/4}}
\frac{dx}{\sqrt{k+\sec(x)^{2}}}=\frac{1}{2\pi}\operatorname{arcsec}\left(
k+1\right)  ,
\]%
\[
T_{k}=\frac{\sqrt{k}}{\pi^{2}}%
{\displaystyle\int\limits_{0}^{\pi/4}}
\,%
{\displaystyle\int\limits_{0}^{\pi/4}}
\frac{dx\,dy}{\sqrt{k+\sec(x)^{2}+\sec(y)^{2}}}=\frac{1}{2\pi^{2}}%
{\displaystyle\int\limits_{0}^{\pi\,S_{k}}}
\operatorname{arcsec}\left(  1+\frac{k(k+1)}{k-\tan(z)^{2}}\right)  dz,
\]%
\[%
\begin{array}
[c]{ccc}%
U=\dfrac{1}{\pi^{2}}%
{\displaystyle\int\limits_{0}^{1}}
\dfrac{\operatorname{arcsec}\left(  2t^{2}+4\right)  }{\left(  2t^{2}%
+1\right)  \sqrt{2t^{2}+3}}dt, &  & V=\dfrac{1}{\pi^{2}}%
{\displaystyle\int\limits_{0}^{1}}
\dfrac{\operatorname{arcsec}\left(  t^{2}+5\right)  }{\left(  t^{2}+2\right)
\sqrt{t^{2}+4}}dt.
\end{array}
\]
The preceding table complements an analogous table in \cite{Fi2} for first and
second moments of $\max\{X_{1},X_{2},...,X_{n}\}$. Similar expressions for
$\mu_{8}=2.847...$ and $\nu_{8}=8.778...$ remain to be found.

\section{Key Observation}

Let us rescale length so that the circumradius of the $n$-simplex is $1$.
Adjusted width will be denoted by $\widetilde{w}_{n}$. Using optimality
properties of the $n$-simplex, Sun \cite{Sn} deduced a formula for mean half
width:%
\begin{align*}
\frac{1}{2}\mathbb{E}\left(  \widetilde{w}_{n}\right)   & =\frac{n+1}{2}%
\sqrt{\frac{(n+1)n}{2\pi}}\frac{\Gamma\left(  \frac{n}{2}\right)  }%
{\Gamma\left(  \frac{n+1}{2}\right)  }%
{\displaystyle\int\limits_{-\infty}^{\infty}}
F\left(  \tfrac{x}{\sqrt{2}}\right)  ^{n-1}f(x)\,dx\\
& =\frac{(n+1)^{3/2}}{\sqrt{2n}}\frac{\Gamma\left(  \frac{n}{2}\right)
}{\Gamma\left(  \frac{n+1}{2}\right)  }%
{\displaystyle\int\limits_{-\infty}^{\infty}}
x\,F(x)^{n}f(x)\,dx
\end{align*}
(see Corollary 2 on p. 1581 and its proof on p. 1585; his $M$ is the same as
our $n+1$). \ We recognize the latter integral as $\mu_{n+1}/(2(n+1))$; hence
\[
\mathbb{E}\left(  \widetilde{w}_{n}\right)  =\sqrt{\frac{n+1}{2n}}\frac
{\Gamma\left(  \frac{n}{2}\right)  }{\Gamma\left(  \frac{n+1}{2}\right)  }%
\mu_{n+1}
\]
and therefore%
\[
\mathbb{E}\left(  w_{n}\right)  =\frac{1}{2}\frac{\Gamma\left(  \frac{n}%
{2}\right)  }{\Gamma\left(  \frac{n+1}{2}\right)  }\mu_{n+1}
\]
because, in our original scaling, the circumradius is $\sqrt{n/(2(n+1))}$.

No similar integral expression for $\mathbb{E}\left(  \widetilde{w}_{n}%
^{2}\right)  $ appears in \cite{Sn}. We circumvent this difficulty by noticing
that the formula \cite{BY, Lc}
\[
\mathbb{E}\left(  \sqrt{%
{\displaystyle\sum\limits_{k=1}^{n}}
X_{k}^{2}}\right)  =\sqrt{2}\frac{\Gamma\left(  \frac{n+1}{2}\right)  }%
{\Gamma\left(  \frac{n}{2}\right)  }
\]
bears some resemblance to the coefficient of $\mu_{n+1}$ in our expression for
$\mathbb{E}\left(  w_{n}\right)  $. The square version%
\[
\mathbb{E}\left(
{\displaystyle\sum\limits_{k=1}^{n}}
X_{k}^{2}\right)  =n
\]
is trivial and leads us to conjecture that
\[
\mathbb{E}\left(  w_{n}^{2}\right)  =\frac{1}{2n}\nu_{n+1}
\]
by analogy. \ Numerical confirmation for $n\leq6$ is possible via the computer
algebra technique described in \cite{Fi1}.

In summary, we have mean width results%
\[
\mathbb{E}\left(  w_{2}\right)  =\frac{3}{\pi}=0.954929658551372...,
\]%
\[
\mathbb{E}\left(  w_{3}\right)  =\frac{3}{2}\left(  1-2S_{2}\right)
=0.912260171954089...,
\]%
\[
\mathbb{E}\left(  w_{4}\right)  =\frac{20}{3\pi}\left(  1-3S_{2}\right)
=0.874843256085440...,
\]%
\[
\mathbb{E}\left(  w_{5}\right)  =\frac{45}{16}\left(  1-4S_{2}+2T_{2}\right)
=0.842274297659162...,
\]%
\[
\mathbb{E}\left(  w_{6}\right)  =\frac{56}{5\pi}\left(  1-5S_{2}%
+5T_{2}\right)  =0.813743951590337...,
\]
and mean square width results%
\[
\mathbb{E}\left(  w_{2}^{2}\right)  =\frac{1}{2}\left(  1+\frac{3\sqrt{3}%
}{2\pi}\right)  =0.913496671566344...,
\]%
\[
\mathbb{E}\left(  w_{3}^{2}\right)  =\frac{1}{3}\left(  1+\frac{3+\sqrt{3}%
}{\pi}\right)  =0.835419517991054...,
\]%
\begin{align*}
\mathbb{E}\left(  w_{4}^{2}\right)   & =\frac{1}{4}\left(  1+\frac{5\sqrt{3}%
}{2\pi}+\frac{30}{\pi}S_{1/2}-\frac{5\sqrt{3}}{\pi}S_{3}\right) \\
& =0.769572883591771...,
\end{align*}%
\begin{align*}
\mathbb{E}\left(  w_{5}^{2}\right)   & =\frac{1}{5}\left(  1+\frac
{5(9+2\sqrt{3})}{2\pi}-\frac{90}{\pi}S_{2}-\frac{15\sqrt{3}}{\pi}S_{3}\right)
\\
& =0.714241915072694...,
\end{align*}%
\begin{align*}
\mathbb{E}\left(  w_{6}^{2}\right)   & =\frac{1}{6}\left(  1+\frac{35\sqrt{3}%
}{4\pi}+\frac{210}{\pi}S_{1/2}-\frac{105}{\pi}S_{2}-\frac{35\sqrt{3}}{\pi
}S_{3}\right. \\
& \;\;\;\;\;\;\;\;\left.  +\frac{35\sqrt{3}}{2\pi}T_{3}+\frac{210}{\pi}%
U-\frac{420}{\pi}V\right) \\
& =0.667314714095430....
\end{align*}

\section{Asymptotics}

We turn now to the asymptotic distribution of $r_{n}$ as $n\rightarrow\infty$.
Define $a_{n}$ to be the positive solution of the equation \cite{Da, Ha}%
\[
2\pi\,a_{n}^{2}\exp\left(  a_{n}^{2}\right)  =n^{2},
\]
that is,%
\[
a_{n}=\sqrt{W\left(  \frac{n^{2}}{2\pi}\right)  }\sim\sqrt{2\ln(n)}-\frac
{1}{2}\frac{\ln(\ln(n))+\ln(4\pi)}{\sqrt{2\ln(n)}}
\]
in terms of the Lambert $W$ function \cite{Cr}. \ It can be proved that the
required density is a convolution \cite{Gb, Cx}:
\begin{align*}
\lim_{n\rightarrow\infty}\frac{d}{dy}\operatorname*{P}\left(  \sqrt{2\ln
(n)}(r_{n}-2a_{n})<y\right)   & =%
{\displaystyle\int\limits_{-\infty}^{\infty}}
\exp(-x-e^{-x})\exp(-(y-x)-e^{-(y-x)})dx\\
& =2\,e^{-y}K_{0}\left(  2\,e^{-y/2}\right)
\end{align*}
where $K_{0}$ is the modified Bessel function of the second kind \cite{Ov}. A
random variable $Y$, distributed as such, satisfies
\[%
\begin{array}
[c]{ccc}%
\mathbb{E}(Y)=2\gamma, &  & \mathbb{E}(Y^{2})=\dfrac{\pi^{2}}{3}+4\gamma^{2}%
\end{array}
\]
where $\gamma$ is the Euler-Mascheroni constant \cite{Fi3}. \ This implies
that
\[
\mu_{n}\sim2\left(  a_{n}+\frac{\gamma}{\sqrt{2\ln(n)}}\right)  \sim
2\sqrt{2\ln(n)}-\frac{\ln(\ln(n))+\ln(4\pi)-2\gamma}{\sqrt{2\ln(n)}}
\]
and hence%
\begin{align*}
\mathbb{E}(w_{n})  & =\frac{1}{2}\frac{\Gamma\left(  \frac{n}{2}\right)
}{\Gamma\left(  \frac{n+1}{2}\right)  }\mu_{n+1}=\frac{1}{2}\frac
{\Gamma\left(  \frac{n}{2}\right)  }{\Gamma\left(  \frac{n+1}{2}\right)
}\cdot\frac{\mu_{n+1}}{\mu_{n}}\cdot\mu_{n}\\
& \sim\frac{1}{\sqrt{2n}}\left(  1+\frac{1}{4n}\right)  \cdot\left(
1+\frac{1}{2n\,\ln(n)}\right)  \cdot2\left(  a_{n}+\frac{\gamma}{\sqrt
{2\ln(n)}}\right) \\
& \sim2\sqrt{\frac{\ln(n)}{n}}-\frac{\ln(\ln(n))+\ln(4\pi)-2\gamma}%
{2\sqrt{n\ln(n)}}.
\end{align*}
More terms in the asymptotic expansion are possible.

If we rescale length so that the inradius of the $n$-simplex is $1$ and denote
adjusted width by $\widehat{w}_{n}$, then%
\[
\mathbb{E}(\widehat{w}_{n})\sim\sqrt{2}n\cdot2\sqrt{\frac{\ln(n)}{n}}%
\sim2\sqrt{2n\ln(n)}
\]
because, in our original scaling, the inradius is $\sqrt{1/(2n(n+1))}$. \ This
first-order approximation is consistent with \cite{BS}.

\section{Regular Octahedron}

As an aside, we return to the setting of $\mathbb{R}^{3}$ and review our
computational methods for $C=$ the regular octahedron with edges of unit length.

For simplicity, let $\Diamond$ be the octahedron with vertices%
\[%
\begin{array}
[c]{lllll}%
v_{1}=\left(  1,0,0\right)  , &  & v_{2}=\left(  -1,0,0\right)  , &  &
v_{3}=\left(  0,1,0\right)  ,\\
v_{4}=\left(  0,-1,0\right)  , &  & v_{5}=\left(  0,0,1\right)  , &  &
v_{6}=\left(  0,0,-1\right)  .
\end{array}
\]
At the end, it will be necessary to normalize by $\sqrt{2}$, the edge-length
of $\Diamond$.

Also let $\widetilde{\Diamond}$ be the union of six overlapping balls of
radius $1/2$ centered at $v_{1}/2$, $v_{2}/2$, $v_{3}/2$, $v_{4}/2$, $v_{5}/2
$, $v_{6}/2$. Clearly $\Diamond\subset\widetilde{\Diamond}$ and
$\widetilde{\Diamond}$ has centroid $(0,0,0)$. \ A diameter of
$\widetilde{\Diamond}$ is the length of the intersection between
$\widetilde{\Diamond}$ and a line passing through the origin.

Computing all widths of $\Diamond$ is equivalent to computing all diameters of
$\widetilde{\Diamond}$. The latter is achieved as follows. Fix a point
$(a,b,c)$ on the unit sphere. \ The line $L$ passing through $(0,0,0)$ and
$(a,b,c)$ has parametric representation%
\[%
\begin{array}
[c]{ccccccc}%
x=t\,a, &  & y=t\,b, &  & z=t\,c, &  & t\in\mathbb{R}%
\end{array}
\]
and hence $y=(b/a)x$, $z=(c/a)x$ assuming $a\neq0$. \ The nontrivial
intersection between first sphere and $L$ satisfies%
\[
\left(  x-\tfrac{1}{2}\right)  ^{2}+\left(  \tfrac{b}{a}x\right)  ^{2}+\left(
\tfrac{c}{a}x\right)  ^{2}=\tfrac{1}{4}
\]
thus $x_{1}=a^{2}$ since $a^{2}+b^{2}+c^{2}=1$; the nontrivial intersection
between second sphere and $L$ satisfies%
\[
\left(  x+\tfrac{1}{2}\right)  ^{2}+\left(  \tfrac{b}{a}x\right)  ^{2}+\left(
\tfrac{c}{a}x\right)  ^{2}=\tfrac{1}{4}
\]
thus $x_{2}=-a^{2}$. \ The nontrivial intersection between third/fourth sphere
and $L$ satisfies%
\[
x^{2}+\left(  \tfrac{b}{a}x\mp\tfrac{1}{2}\right)  ^{2}+\left(  \tfrac{c}%
{a}x\right)  ^{2}=\tfrac{1}{4}
\]
thus $x_{3}=a\,b$ , $x_{4}=-a\,b$. The nontrivial intersection between
fifth/sixth sphere and $L$ satisfies%
\[
x^{2}+\left(  \tfrac{b}{a}x\right)  ^{2}+\left(  \tfrac{c}{a}x\mp\tfrac{1}%
{2}\right)  ^{2}=\tfrac{1}{4}
\]
thus $x_{5}=a\,c$ , $x_{6}=-a\,c$. \ \ \ 

We now examine all pairwise distances, squared, between the six intersection
points:%
\begin{align*}
& \left(  x_{i}-x_{j}\right)  ^{2}+\left(  \tfrac{b}{a}x_{i}-\tfrac{b}{a}%
x_{j}\right)  ^{2}+\left(  \tfrac{c}{a}x_{i}-\tfrac{c}{a}x_{j}\right)  ^{2}\\
& =\left\{
\begin{array}
[c]{lll}%
4a^{2} &  & \text{if }i=1,j=2\\
1-2a\,b-c^{2} &  & \text{if }i=1,j=3\text{ or }i=2,j=4\text{ }\\
1+2a\,b-c^{2} &  & \text{if }i=1,j=4\text{ or }i=2,j=3\\
1-2a\,c-b^{2} &  & \text{if }i=1,j=5\text{ or }i=2,j=6\\
1+2a\,c-b^{2} &  & \text{if }i=1,j=6\text{ or }i=2,j=5\\
4b^{2} &  & \text{if }i=3,j=4\\
(b-c)^{2} &  & \text{if }i=3,j=5\text{ or }i=4,j=6\\
(b+c)^{2} &  & \text{if }i=3,j=6\text{ or }i=4,j=5\\
4c^{2} &  & \text{if }i=5,j=6
\end{array}
\right.
\end{align*}
and define%
\begin{align*}
g(a,b)  & =\max\left\{  4a^{2},1-2a\,b-c^{2},1+2a\,b-c^{2},1-2a\,c-b^{2}%
,1+2a\,c-b^{2}\right. \\
& \left.  \;\;\;\;\;\;\;\;\;4b^{2},(b-c)^{2},(b+c)^{2},4c^{2}\right\}  .
\end{align*}
The mean width for $C$ is%
\[
\frac{1}{\sqrt{2}}\frac{1}{4\pi}%
{\displaystyle\int\limits_{0}^{2\pi}}
{\displaystyle\int\limits_{0}^{\pi}}
\sqrt{g(\cos\theta\sin\varphi,\sin\theta\sin\varphi,\cos\varphi)}\sin
\varphi\,d\varphi\,d\theta=\frac{3}{\pi}\arccos\left(  \frac{1}{3}\right)
\]
and the mean square width is%
\[
\frac{1}{2}\frac{1}{4\pi}%
{\displaystyle\int\limits_{0}^{2\pi}}
{\displaystyle\int\limits_{0}^{\pi}}
g(\cos\theta\sin\varphi,\sin\theta\sin\varphi,\cos\varphi)\sin\varphi
\,d\varphi\,d\theta=\frac{2}{3}\left(  1+\frac{2\sqrt{3}}{\pi}\right)  .
\]

Here are details on the final integral. A plot of the surface
\[
(\theta,\varphi)\longmapsto\sqrt{\frac{g(\cos\theta\sin\varphi,\sin\theta
\sin\varphi,\cos\varphi)}{2}}
\]
appears in Figure 1, where $0\leq\theta\leq2\pi$ and $0\leq\varphi\leq\pi$.
\ Figure 2 contains the same surface, but viewed from above. \ Our focus will
be on the part of the surface to the right of the bottom center, specifically
$0\leq\theta\leq\pi/4$ and $\pi/2\leq\varphi\lessapprox9/4$. \ The volume
under this part is $1/24^{\text{th}}$ of the volume under the full surface.

We need to find the precise upper bound on $\varphi$ as a function of $\theta
$. Recall the formula for\ $g$ as a maximum over nine terms; let $g_{\ell}$
denote the $\ell^{\text{th}}$ term, where $1\leq\ell\leq9$. \ Then the upper
bound on $\varphi$ is found by solving the equation%
\[
g_{1}(\cos\theta\sin\varphi,\sin\theta\sin\varphi,\cos\varphi)=g_{9}%
(\cos\theta\sin\varphi,\sin\theta\sin\varphi,\cos\varphi)
\]
for $\varphi$. We obtain $\varphi(\theta)=2\arctan(h(\theta))$, where
\[
h(\theta)=\cos\theta+\sqrt{\frac{3+\cos(2\theta)}{2}}
\]
and, in particular,
\[
\varphi(0)=2\arctan\left(  1+\sqrt{2}\right)  \approx2.3562,
\]%
\[
\varphi(\pi/4)=2\arctan\left(  \left(  1+\sqrt{3}\right)  /\sqrt{2}\right)
\approx2.1862.
\]
It follows that $g=g_{1}$ for $0\leq\theta\leq\pi/4$ and $\pi/2\leq\varphi
\leq2\arctan(h)$. \ Now we have%
\begin{align*}
& \frac{1}{2}\frac{1}{4\pi}%
{\displaystyle\int}
g_{1}(\cos\theta\sin\varphi,\sin\theta\sin\varphi,\cos\varphi)\sin
\varphi\,d\varphi\\
& =\frac{\left(  1+\cos(2\theta)\right)  \left(  \cos(3\varphi)-9\cos
(\varphi)\right)  }{48\pi}%
\end{align*}
and
\[
\left.  \cos(3\varphi)\right\vert _{\pi/2}^{2\arctan(h)}=\frac{\left(
1+4h+h^{2}\right)  \left(  1-4h+h^{2}\right)  \left(  1-h^{2}\right)
}{(1+h^{2})^{3}},
\]%
\[
\left.  \cos(\varphi)\right\vert _{\pi/2}^{2\arctan(h)}=\dfrac{1-h^{2}%
}{1+h^{2}},
\]
therefore
\begin{align*}
& \frac{1}{2}\frac{1}{4\pi}%
{\displaystyle\int_{\pi/2}^{2\arctan(h)}}
g_{1}(\cos\theta\sin\varphi,\sin\theta\sin\varphi,\cos\varphi)\sin
\varphi\,d\varphi\\
& =\frac{\left(  h^{4}+4h^{2}+1\right)  \left(  h^{2}-1\right)  \left(
1+\cos(2\theta)\right)  }{6\pi(1+h^{2})^{3}}.
\end{align*}
Integrating this expression from $0$ to $\pi/4$ gives the desired formula for
$\operatorname*{E}\left(  w_{\text{octa}}^{2}\right)  $.

\section{$n$-Cubes}

After having written the preceding, we discovered \cite{HZ}, which gives the
mean width for a regular $n$-simplex in $\mathbb{R}^{n}$ as
\[
\mathbb{E}\left(  w_{n}\right)  =\frac{n(n+1)}{\sqrt{2}\pi}\frac{\Gamma\left(
\frac{n}{2}\right)  }{\Gamma\left(  \frac{n+1}{2}\right)  }%
{\displaystyle\int\limits_{-\infty}^{\infty}}
e^{-2x^{2}}\left(  \frac{1+\operatorname{erf}(x)}{2}\right)  ^{n-1}\,dx.
\]
Consistency is readily established; nothing is said in \cite{HZ} about the
connection between $\mathbb{E}\left(  w_{n}\right)  $ and order statistics
from a normal distribution (more precisely, the expected range $\mu_{n+1}$).

By contrast, the mean width for an $n$-cube with edges of unit length is
elementary:
\[
\mathbb{E}\left(  w_{n\text{-cube}}\right)  =\frac{n}{\sqrt{\pi}}\frac
{\Gamma\left(  \frac{n}{2}\right)  }{\Gamma\left(  \frac{n+1}{2}\right)  }
\]
and we conjecture that%
\[
\mathbb{E}\left(  w_{n\text{-cube}}^{2}\right)  =1+\frac{2(n-1)}{\pi}.
\]

\section{$n$-Crosspolytopes}

A\ regular $n$-crosspolytope with edges of unit length has mean width
\cite{HZ, BH}
\[
\mathbb{E}\left(  w_{n\text{-crosspolytope}}\right)  =\frac{2\sqrt{2}%
n(n-1)}{\pi}\frac{\Gamma\left(  \frac{n}{2}\right)  }{\Gamma\left(  \frac
{n+1}{2}\right)  }%
{\displaystyle\int\limits_{0}^{\infty}}
e^{-2x^{2}}\operatorname{erf}(x)^{n-2}\,dx;
\]
the case $n=3$ corresponds to the octahedron discussed earlier. \ It is not
surprising that a connection exists with order statistics from a half-normal
(folded) distribution. \ We will examine this later, as well as relevant
expressions from \cite{Gv}.\ An appropriate mean square conjecture also needs
to be formulated in this scenario.

\bigskip%
\begin{figure}[ptb]%
\centering
\includegraphics[
height=4.2713in,
width=5.348in
]%
{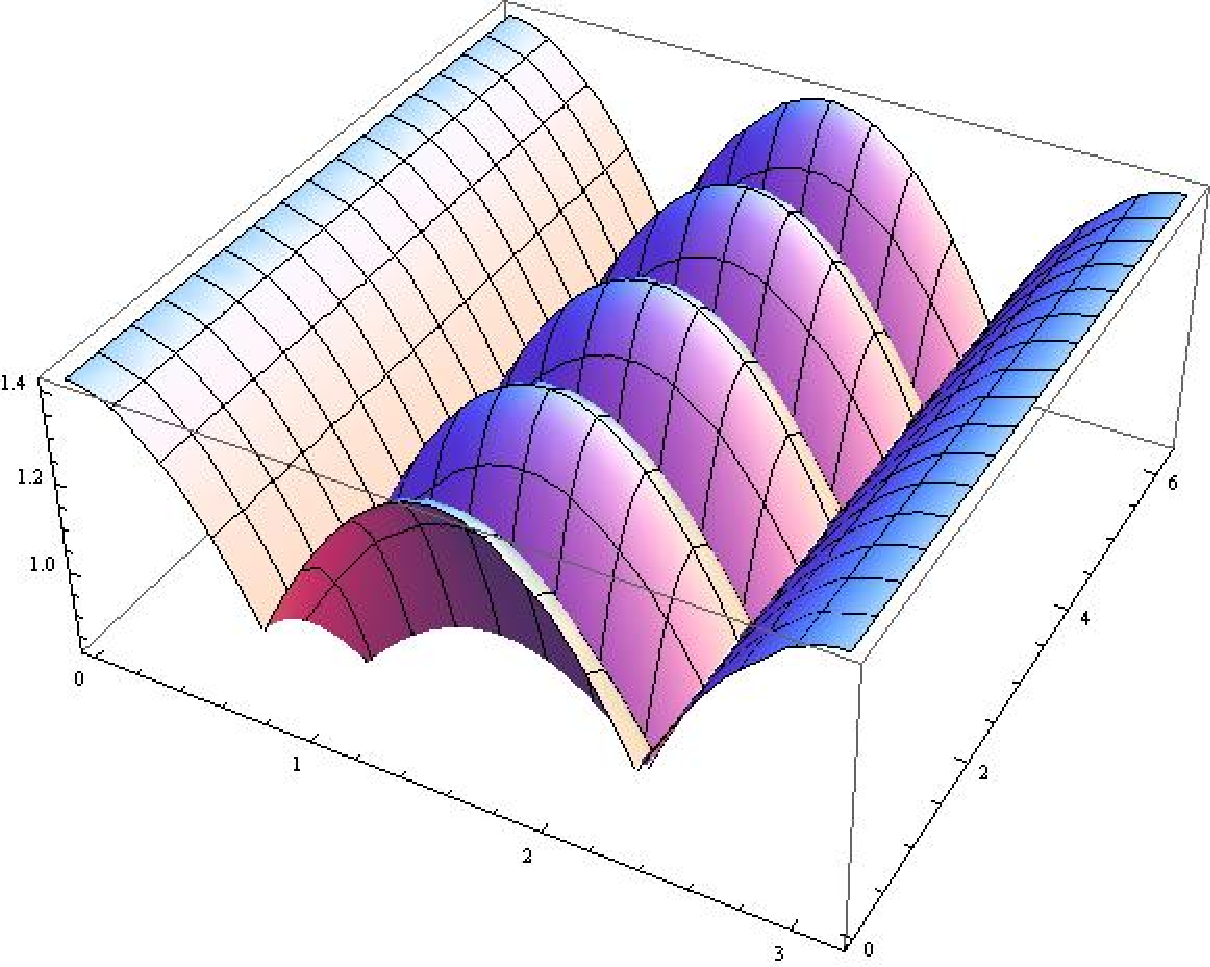}%
\caption{Surface plot of $\sqrt{g/2}$, where $0\leq\theta\leq2\pi$ and
$0\leq\varphi\leq\pi$.}%
\end{figure}
%

\begin{figure}[ptb]%
\centering
\includegraphics[
height=4.2704in,
width=5.2935in
]%
{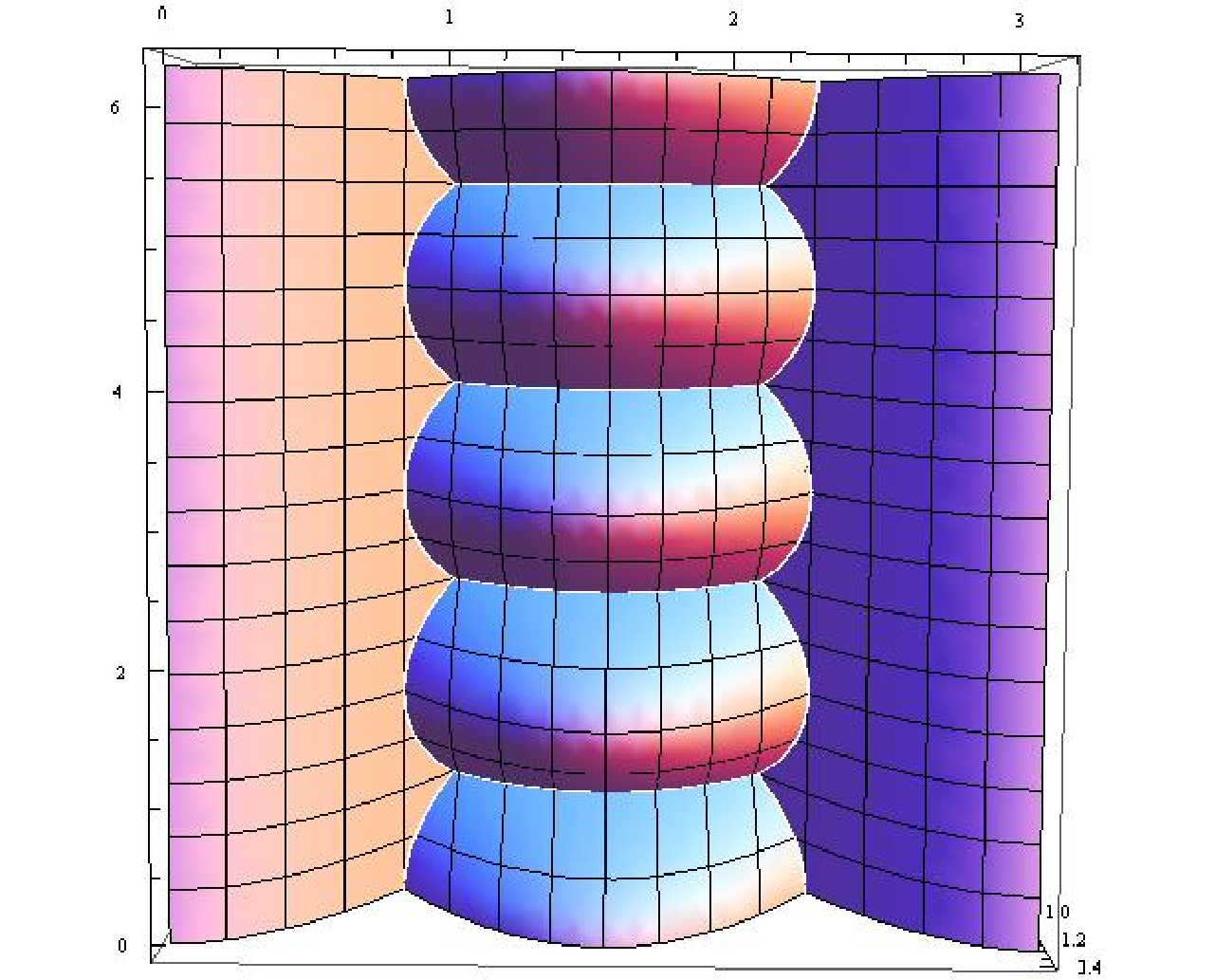}%
\caption{Another view of $\sqrt{g/2}$, with contours of intersection.}%
\end{figure}

\bigskip
\end{document}